\documentclass[10pt]{amsart}
\input epsf
\usepackage{amsmath}
\usepackage{amssymb}
\usepackage{amscd}
\usepackage{amsthm}
\usepackage{subfigure}

\usepackage[all]{xy}
\usepackage{color}

\numberwithin{equation}{section}

\newtheorem{theorem}[equation]{Theorem}

\newtheorem{proposition}[equation]{Proposition}
\newtheorem{prop}[equation]{Proposition}

\newtheorem{lem}[equation]{Lemma}
\newtheorem{lemma}[equation]{Lemma}

\newtheorem{corollary}[equation]{Corollary}

\newtheorem{conjecture}[equation]{Conjecture}

\theoremstyle{remark}
\newtheorem{remark}[equation]{Remark}

\newtheorem{rem}[equation]{Remark}

\theoremstyle{definition}
\newtheorem{definition}[equation]{Definition}

\def\XXint#1#2#3{{\setbox0=\hbox{$#1{#2#3}{\int}$} 
	\vcenter{\hbox{$#2#3$}}\kern-.5\wd0}}

\newcommand{\To}{\longrightarrow}

\newcommand{\N}{\mathbb N}

\newcommand{\R}{\mathbb R}

\newcommand{\diam}{\operatorname{diam}}

\newcommand{\norm}[1]{\left\Vert#1\right\Vert}

\def\eps{\epsilon}

\begin{document}

\title{Lipschitz and path isometric embeddings\\ of metric spaces}

\author{Enrico Le Donne}


\renewcommand{\subjclassname}{%
\textup{2010} Mathematics Subject Classification}

\date{August 28, 2012}

\begin{abstract}
We prove that each sub-Riemannian manifold  can be embedded in some Euclidean space preserving the length of all the curves in the manifold. The result is an extension of Nash $C^1$ Embedding Theorem. 
For more general metric spaces the same result is false, e.g., for Finsler non-Riemannian manifolds. However, we also show that any metric space of finite packing dimension can be embedded in some Euclidean space via a Lipschitz map.
\end{abstract}

\maketitle

\section{Overview}
A map $f:X\to Y$ between two metric spaces 
$X$ and $Y$ is called a {\em path isometry} (probably a better name is a {\em length-preserving map}) if, for all curves $\gamma$ in $ X$, one has
$$L_Y(f\circ\gamma)=L_X(\gamma).$$ 
Here $L_X$ and $L_Y$ denote the lengths of the parametrized curves with respect to the distances of $X$ and of $Y$, respectively.
From the definition, a path isometry is not necessarily injective.

The first aim of the following paper is to show that any sub-Riemannian manifold can be mapped into some Euclidean space via a path isometric embedding, i.e., a topological embedding that is also a path isometry. Sub-Riemannian manifolds are  metric spaces when endowed with the Carnot-Carath\'eodory distance $d_{CC}$ associated to the fixed sub-bundle and Riemannian structure. For an introduction to sub-Riemannian geometry see \cite{bellaiche, Gromov,  Burago,  Montgomery, Buliga1,LeDonne-lectures}.

An equivalent statement of our first result is the following. Denote by $\mathbb E^k$  the $k$-dimensional Euclidean space. Our result says that, for every  sub-Riemannian manifold $(M,d_{CC})$, 
there exists a path connected subset $\Sigma\subset\mathbb E^k$, for some $k\in\N$, such that, when $\Sigma$ is endowed with the path distance $d_\Sigma$ induced by the Euclidean length, then the metric space $(\Sigma, d_\Sigma)$ is isometric to $(M,d_{CC})$.

After such a fact one should wonder which are the length metric spaces obtained as subsets of $\mathbb E^k$ with induced length structure.
We show that any distance on $\R^n$ that comes	 from a norm but not from a scalar product cannot be obtained in such a way.

We conclude the paper by showing another positive result for general metric spaces: every metric space of finite packing dimension 
has a Lipschitz embedding into some $\mathbb E^k$.

\section{Old and new results}

In 1954 John Nash showed that any Riemannian manifold can be seen as a $C^1$ submanifold of some Euclidean space. Namely, for any $n$-dimensional Riemannian manifold
$(M,g)$, there exists a $C^1$ submanifold $N$ of the $(2 n +1)$-dimensional  Euclidean space $\mathbb E^{2 n +1}$ such that $N$, endowed with the restriction of the Euclidean Riemannian tensor, is $C^1$ equivalent to $(M,g)$.
Two Riemannian manifolds $(M_1,g_1)$ and $(M_2,g_2)$ are considered $C^1$ equivalent if there exists a $C^1$ diffeomorphism $f:M_1\to M_2$ such that the pull-back tensor $f^*g_2$ equals $g_1$.
In Riemannian geometry, a $C^1$ map $f$ between two Riemannian manifolds $(M_1,g_1)$ and $(M_2,g_2)$ with the property that 
$$f:(M_1,g_1)\to(f(M_1),g_2|_{T(f(M_1))})$$
is a $C^1$ equivalence is called an `isometric embedding'.
However, in the present paper we will avoid such a term for the reason that  the notion of isometric embedding is different in the setting of metric spaces. Indeed, let $d_{g_1}$ and $d_{g_2}$ be the distance functions on $M_1$ and $M_2$, respectively, induced by $g_1$ and $g_2$, respectively.
Then the fact that $f:(M_1,g_1)\to(M_2,g_2)$ is a Riemannian `isometric embedding'   does not imply that 
$f:(M_1,d_{g_1})\to(M_2,d_{g_2})$ is an isometric embedding of the metric space $(M_1,d_{g_1})$ into the metric space $(M_2,d_{g_2})$, i.e., it is not true in general that
$$ d_{g_2}(f(p),f(q))=d_{g_1}(p,q),\qquad \forall p,q\in M_1.$$
However, an elementary but important consequence of having a Riemannian isometric embedding is that the length of paths is preserved. In other words, Nash's theorem can be restated as saying that any Riemannian manifold can be path isometrically embedded into some Euclidean space.

%
%
\begin{definition}[Path isometric embedding]
A map $f:X\to Y$ between two metric spaces 
$X$ and $Y$ is called a {\em path isometric embedding} if it is a topological embedding, i.e., a homeomorphism onto its image, and, for all curves $\gamma\subset X$, one has
$$L_Y(f\circ \gamma)=L_X(\gamma).$$ 
\end{definition}
We want to clarify that the above condition is required also for curves of infinite length.

One of the versions of Nash Theorem can be stated as follows.
\begin{theorem}[Nash]\label{Nash0}
Let $(M,g)$ be a $C^\infty$ Riemannian manifold of dimension $n$.
Then there exists a $C^1$ path isometric embedding
$$f:(M,d_g)\to \mathbb E^k,$$
with $k=2 n+1$.
\end{theorem}

The theorem  originally appeared in \cite{Nash}, later it was generalized by Nicolaas Kuiper in \cite{Kuiper}.
The Nash-Kuiper $C^1$ Theorem can be stated  in the following form.

\begin{theorem}[Nash-Kuiper $C^1$ Embedding Theorem]\label{Nash-Kuiper}
Let $(M,g)$ be a $C^\infty$ Riemannian manifold of dimension $n$.
If there is 
a $C^\infty$ $1$-Lipschitz embedding
$$f:(M,d_g)\to \mathbb E^k$$
into an Euclidean space $\mathbb E^k$ with $k\geq n+1$,
 then, for all $\eps>0$, there exists a 
$C^1$ path isometric embedding
$$\bar f:(M,d_g)\to \mathbb E^k,$$
that is $\eps$-close to $f$, i.e., for any $p\in M$,
$$d_{\mathbb E}(f(p), \bar f(p))\leq \eps.$$
\end{theorem}

In particular, as follows from a result of Nash which extends the Whitney Embedding Theorem, any $n$-dimensional Riemannian manifold admits a path isometric $C^1$ embedding into an arbitrarily small neighborhood in $(2 n+1)$-dimensional Euclidean space. 

The Nash-Kuiper Theorem has many counter-intuitive implications. For example, 
it follows that  there exist $C^1$ path isometric embeddings of the hyperbolic plane in $\mathbb E^3$. Additionally, 
any closed, oriented Riemannian surface can be $C^1$ path isometrically embedded into an arbitrarily small ball in $\mathbb E^3$. 
 Whereas, for curvature reasons, there is no such a $C^2$-embedding.

In \cite[2.4.11]{GromovPDR} Gromov proved that any Riemannian manifold of dimension $n$ admits a path isometry  into $\mathbb E^n$ (notice the same dimension).
In a recent paper  \cite{Petrunin}  Petrunin extended Gromov's result to sub-Riemannian manifolds for a more rigid class of maps: the intrinsic isometries.
The key fact used by Petrunin is that any sub-Riemannian distance is a {\em monotone} limit of Riemannian distances. 
Such a fact is well known in nonholonomic geometry since the last 25 years, and probably is due to V. Gershkovich. 
 This observation will be essential in considering limits of Nash's embeddings as we will do in this paper.

For topological reasons, both Gromov's and Petrunin's maps are in general not injective. Our aim is to have path isometries that are also embeddings. 
Nontheless, this paper has been strongly influenced by the work of Petrunin. Some of the methods are just elaborations and generalizations of Petrunin's ideas.
As an example of the fact that Petrunin's notion of intrinsic isometry is related with our work, we shall show that 
 any path isometric embedding is an intrinsic isometry, cf.~Section \ref{intrinsic isometry section}.
 
As a first result, we provide a generalization of Nash Theorem
 to metric spaces obtained as limit of an increasing sequence of Riemannian metrics on a fixed  manifold, e.g., sub-Riemannian manifolds.
\begin{theorem}[Path Isometric Embedding] \label{path-isom-embed}
Let $M$ be a $C^\infty$  manifold of dimension $n$.
Let $g_m$ be a sequence of Riemannian structures on $M$ and let $d_{g_m}$ be the distance function induced by $g_m$.
Assume that, for all $p$ and $q\in M$, for all $m\in \N$,
$$d_{g_m}(p,q)\leq d_{g_{m+1}}(p,q).$$
Assume also that, for all $p$ and $q\in M$, the limit
$$d(p,q):=\lim_{m\to\infty}d_{g_m}(p,q)$$
is finite 
 and that the function $d$ gives a distance that induces the manifold topology on $M$.
Then there exists a   path isometric embedding
$$f:(M,d)\to \mathbb E^k,$$
with $k=2 n+1$.
\end{theorem}

In Section \ref{subRiem Sect} we will recall the general  definition of a sub-Riemannian manifold
 and show that  the sub-Riemannian distance function is a point-wise  limit of Riemannian distance functions.
Then  the following fact will be an immediate consequence of the above theorem.
\begin{corollary}\label{subRiemEmb}
Each sub-Riemannian manifold of topological dimension $n$ can be path isometrically embedded into $\mathbb E^{2 n+1}$.
\end{corollary}

Actually, the proof of Theorem \ref{path-isom-embed} gives a more precise result for the dimension of the target. 
\begin{corollary}\label{subRiemEmb2}
As in Theorem \ref{path-isom-embed}, let $(M,d)$ be a metric space obtained as a 
limit of an increasing sequence of Riemannian metrics on a    manifold of topological dimension $n$.
Let $d_{\rm Riem}$ be some Riemannian distance such that
$$d_{\rm Riem}\leq d.$$
If there exists a $C^\infty$ $1$-Lipschitz embedding
$$f:(M, d_{\rm Riem})\to \mathbb E^k$$
into an Euclidean space $\mathbb E^k$ with $k\geq n+1$,
 then there exists a 
 path isometric embedding
$$\bar f:(M,d)\to \mathbb E^k.$$
\end{corollary}
Consequently, the Heisenberg group endowed with the usual Carnot-Carath\'eodory metric is isometric to a subset of $\R^4$ endowed with the path metric induced by the Euclidean distance, cf.~Corollary \ref{HinR4}.
Similarly, the Grushin plane  can be realized as a subset of $\R^3$ with the induced path distance.

Our result does not contradict the
 biLipschitz non-embeddability of Carnot-Carath\'eodory spaces.
 Let us recall that it was observed by Semmes 
\cite[Theorem 7.1]{Semmes2} that
Pansu's version of Rademacher's Differentiation Theorem  \cite{Pansu,Margulis-Mostow} implies that
 a Lipschitz embedding of a manifold $M$ endowed with a  sub-Riemannian distance induced by a regular distribution
 into an Euclidean space
cannot be biLipschitz, unless $M$ is in fact Riemannian.
Indeed, in the case of the Heisenberg group $\mathbb H$, any Lipschitz map collapses in the
direction of the center, i.e., 
\begin{equation}
\lim_{g\to e}\,\frac{\|f(gx)-f(x)\|_{\mathbb E}}{d_{CC}(gx,x)}= 0\, ,\quad
\quad g\in{\rm Center}(\mathbb H)\, .
\end{equation}
From this fact we understand that any path isometric embedding $f:\mathbb H\to \mathbb E^k$, which is always a Lipschitz map, has the property that,
for $x\in \mathbb H$, as $g$ goes to the identity element inside ${\rm Center}(\mathbb H)$, the point  
$f(gx)$ converges to $f(x)$ in $\mathbb E^k$ faster than $gx$ converges to $x$ in $\mathbb H$. This last fact does not contradict the existence of a curve $\gamma$ inside $f(\mathbb H)$ from $f(gx)$  to $f(x)$ of length exactly $d_{CC}(gx,x)$ and the fact that all the other  curves
inside $f(\mathbb H)$ from $f(gx)$  to $f(x)$ are not shorter, as the path isometric embedding property would imply.

Also, Corollary \ref{subRiemEmb} does not give any dimensional contradiction. Indeed,   the path metric $d_\Sigma$ on a subset $\Sigma\subset\mathbb E^k$ is larger than the restriction on $\Sigma$ of the Euclidean distance. Thus the metric space $(\Sigma, d_\Sigma)$ can a priori have Hausdorff dimension strictly greater than $k=\dim_H(\mathbb E^k)$. The embeddings of Corollary \ref{subRiemEmb} give non-constructive examples of sets $\Sigma\subset\R^k$ with the property that $$\dim_H(\Sigma, d_\Sigma)>k.$$
Notice that for such examples, the metric $ d_\Sigma$ induces on $\Sigma$ the subspace topology    of $\R^k$.
\\

For the sake of completeness let us  mention the following different generalization by   D'Ambra 
 of Nash's result to the case of contact manifolds. 
  Namely,
let $(M_1,\xi_1,g_1)$ and $(M_2,\xi_2,g_2)$ be two contact manifolds with contact structures $\xi_1$ and $\xi_2$, respectively, and Riemannian metrics 
$g_1$ and $g_2$, respectively. The main result in \cite{DAmbra} claims that if $\dim (M_2) \geq 2 \dim (M_1) + 3$ and $M_1$ is compact,
then there exists a $C^1$ embedding $$f:M_1 \to   M_2,$$
preserving the contact structures and the Riemannian tensors on  $\xi_1$, i.e.,
$$ f_*\xi_1\subset\xi_2 \quad\text{ and } g_1|_{\xi_1}=f^*(g_2|_{f_*\xi_1} )  .$$
\\

We consider now possible generalizations of Theorem \ref{path-isom-embed}.
It is not true that any finite dimensional metric space admits a path isometric embedding into some Euclidean space. Indeed, there is no path isometry from $(\R^2, \norm{\cdot}_\infty)$ to any $\mathbb E^k$. Here  $\norm{\cdot}_\infty$ is the supremum norm on $\R^2$, which does not come from a scalar product. 
Such a nonexistence has been previously pointed out for non-Euclidean normed spaces in \cite{Petrunin}. We provide the following generalization.
\begin{proposition} \label{Finsler}
Let $(M,\norm{\cdot})$ be a Finsler manifold.
If there exists a   path isometry
$$f:(M,\norm{\cdot})\to \mathbb E^k,$$ 
then the manifold is in fact Riemannian.
\end{proposition}
The proof of the above proposition is a consequence of Rademacher's Theorem.
A similar argument is in \cite [Proposition 1.7]{Petrunin}.  
We shall give a more general proof in details.
 \\

An important topological theorem, due to K. Menger and G. N\`obeling, states that any compact metrizable space of topological dimension $m$ can be embedded in $\R^k$ for $k = 2m + 1$. For a reference, see \cite{Munkres}. 
 We shall show the analogue  for Lipschitz embeddings of metric spaces, whose proof is an application of the Baire Category Theorem as well as  the topological version of the theorem.
\begin{theorem}[Lipschitz Embedding] \label{Lip-embed}
Any compact metric space of packing dimension $k$ can be embedded in $\mathbb E^N$ via a Lipschitz map, for $N = 2 k + 1$. \end{theorem}

Since  compact sub-Finsler manifolds are biLipschitz equivalent to  sub-Riemannian manifolds,  any  sub-Finsler manifold is  locally  biLipschitz equivalent to  a  subset of some $\mathbb E^k$ with the path distance. In other words, any sub-Finsler manifold can be embedded into $\mathbb E^k$  via a map that distorts lengths by a controlled ratio. Namely, we already know that for sub-Finsler manifolds the following conjecture holds.
Before stating the conjecture, let us recall the definition of bounded-length-distortion maps.
\begin{definition}[BLD] A map $f:X\to Y$ between two metric spaces 
$X$ and $Y$ is said of  {\em bounded-length-distortion} (BLD for short), if there exists a constant $C$ such that, for all curves $\gamma\subset X$, one has
\begin{equation}\label{BLD}
C^{-1}L_X(\gamma)\leq L_Y(f\circ\gamma)\leq C L_X(\gamma).
 \end{equation}
\end{definition}
\begin{conjecture}[BLD embeddings]
Any compact length metric space  of finite Hausdorff dimension   can be embedded in some Euclidean space via a bounded-length-distortion map.
\end{conjecture}
We expect  the above conjecture   to hold, more because of lack of counterexamples than for  actual reasoning. The map given by Theorem \ref{Lip-embed} satisfies the upper bound of equation \eqref{BLD}. However, even if such a map is injective, it might not satisfy the lower bound of equation \eqref{BLD}.

\subsection{Organization of the paper}
In Section \ref{Existence}, after some preliminary results, we give the proof of   Theorem \ref{path-isom-embed}.
Namely, we show the existence of path isometric embeddings for  
metric spaces obtained as limit of an increasing sequence of Riemannian metrics on a fixed  manifold.

In Section \ref{More}, we present the proof of the corollaries of Theorem \ref{path-isom-embed} and some other consequences. 
Namely, we start by recalling the most general definition of sub-Riemannian distances. Then we show that each such a distance can be obtained as limit of an increasing sequence of Riemannian metrics, proving Corollary  \ref{subRiemEmb}.
Then we prove Corollary \ref{subRiemEmb2}, the  more general version of Theorem \ref{path-isom-embed}.
In Proposition \ref{path-isom-isom}, 
 we show that a map is 
 a path isometric embeddings if and only if it is  an isometry   when
 one gives  the image the path metric induced by the ambient space. 
In connection with the work of Petrunin,  in Proposition \ref{intrinsic prop} we show that a path isometric embedding  between proper geodesic spaces is always an intrinsic isometry. We conclude Section \ref{More} by showing the proof of Proposition \ref{Finsler}, i.e., a Finsler manifold cannot be path isometrically embedded in any Euclidean space, unless it is Riemannian.

Section \ref{Lipschitz} is devoted to the  proof of the Embedding Theorem \ref{Lip-embed}. Namely, any   metric space with finite packing dimension can be Lipschitz embedded in some Euclidean space.

\section{Existence of path isometric embeddings}\label{Existence}

\subsection{Preliminaries}
The following Theorem \ref{Nash1} might seem an easy corollary of Nash-Kuiper Theorem \ref{Nash-Kuiper}. Indeed, by Nash-Kuiper,
any smooth $1$-Lipschitz embedding is arbitrarily close to a $C^1$  length-preserving embedding.
By smoothing one can show the following result: any smooth $1$-Lipschitz embedding is arbitrarily close to a $C^\infty$  embedding that distorts lengths by a factor that is arbitrarily close to $1$. However, the claim of  Theorem \ref{Nash1} is one of the strategic steps of Nash-Kuiper's proof, see \cite[Equation 26, page 390]{Nash} and \cite{Kuiper}.
\begin{theorem}[Consequence of Nash's proof]\label{Nash1}
Let $(M,g)$ be a $C^\infty$ Riemannian manifold of dimension $n$.
If there is a $C^\infty$ $1$-Lipschitz embedding
$$f:(M,d_g)\to \mathbb E^k$$
into an Euclidean space $\mathbb E^k$ with $k\geq n+1$, then, for any $a>0$ and for any  continuous function  $b:M\to\R_{>0}$, there exists a 
$C^\infty$ $1$-Lipschitz embedding
$$\bar f:(M,d_g)\to \mathbb E^k,$$
such that, for any curve $\gamma\subset M$,
$$(1-a) L_g(\gamma)\leq L_{\mathbb E}(\bar f\circ\gamma)\leq L_g(\gamma)$$
and, for any $p\in M$,
$$d_{\mathbb E}(f(p), \bar f(p))\leq b(p).$$
\end{theorem}

For compact manifolds the following result is an easy consequence of Whitney Embedding Theorem, where in fact one can take $k=2 n$.
For general manifolds a proof can be found in \cite[page 394]{Nash}.
\begin{theorem}[Whitney-Nash]\label{Nash2}
Let $(M,g)$ be a $C^\infty$ Riemannian manifold of dimension $n$.
Then there exists a $C^\infty$ $1$-Lipschitz embedding
$$f:(M,d_g)\to \mathbb E^k,$$
with $k=2 n+1$.
\end{theorem}

	 Given a set $\Sigma\subset \mathbb E^k$, one can consider the path metric on $\Sigma$ induced by $L_{\mathbb E}$, i.e., for $p,q\in\Sigma$,
define
$$d_\Sigma(p,q):=\inf\left\{L_{\mathbb E}(\gamma)\;|\;{\rm Im}(\gamma)\subset \Sigma,\; \gamma \text{ from } p \text{ to }q\right\}.$$

\begin{remark}
The function $d_\Sigma$ is a distance whose induced topology, a priori, might be different from the topology of $\Sigma$ as subset of  $ \mathbb E^k$. However, the length structures $L_{\mathbb E}$ and $L_{d_\Sigma}$ coincide. Namely, if $\gamma:I\to  (\Sigma,d_\Sigma) $ is a curve   
then
$$ L_{\mathbb E}(\gamma)  =  L_{d_\Sigma}(\gamma).$$
Such an equality is easy to show. A detailed and more general proof can be found in \cite[Proposition 2.3.12]{Burago}.
\end{remark}

The following fact is the key   for preventing loss of length in the limit process while proving Theorem \ref{path-isom-embed}. A similar argument was used in \cite{Petrunin}. 
\begin{definition}[Neighborhood  $I(\delta)$]
 Let $f:M\to \R^k$
be a $C^\infty$  embedding.
Let $\delta:M\to\R_{>0}$ be a continuous function. We consider the $\delta$-neighborhood of $f(M)$ as the set
$$I(\delta):=I_\delta(f(M)):=\{x\in\R^k:\norm{x-f(p)}_{\mathbb E}<\delta(p),\text{ for some }p\}.$$
\end{definition}

\begin{lemma}[Control on tubular neighborhoods]\label{tubular}
Let $M$ be a $C^\infty$  manifold.
Let
$$f:M\to \R^k$$
be a $C^\infty$  embedding.
Then, for any $ a>0$, there exists a positive continuous function $\delta=\delta_{f, a}:M\to(0, a)$   such that, for all $x,y\in   f(M) $,
$$(1- a)d_{f(M)}(x,y)\leq d_{I(\delta)}(x,y)\leq d_{f(M)}(x,y),$$
where   $d_{f(M)}$ and $d_{I(\delta)}$ are the path metrics in $f(M)$ and ${I(\delta)}$, respectively.
\end{lemma}

 Lemma \ref{tubular} is well-known.  One can give an easy    proof   using the  Neighborhood Theorem.
 A reference for the proof is \cite{Federer-Curvature-measures}.

\subsection{Proof of the existence of path isometric embeddings}
This section is devoted to the proof of Theorem \ref{path-isom-embed}. We will first construct the map $f$, then prove that it is a path isometry, and finally that it is an embedding.
\subsubsection{The construction of $f$}
The map $f$ shall be obtained as a limit of maps $f_m$. The construction of the sequence $f_m$ is by induction. Briefly speaking, we have that $f_m$ is an  isometric  embedding  for the Riemannian structure of $g_m$ obtained by $f_{m-1}$, via Nash-Kuiper Embedding Theorem \ref{Nash-Kuiper}, inside a suitably controlled neighborhood.

From Theorem \ref{Nash2}, we can start with a $C^\infty$ $1$-Lipschitz embedding
$$f_1:(M,g_1)\to \mathbb E^k.$$
For $m\in \N$, set 
$$  \quad a_m:=\dfrac{1}{m}  . $$ 
Considering the function $\delta_{f, a}$ of Lemma \ref{tubular}, set $\delta_1:=\delta_{f_1, a_1}$.
Choose any $C^0$ function $b_1$ with $0<b_1(p)<  \delta_1(p)$, for all $p\in M$.

By recurrence, for each $m\in \N$, perform the following construction of 
$C^\infty$ $1$-Lipschitz embeddings
$$f_m:(M,g_m)\to \mathbb E^k$$
and positive  continuous function $b_m$ and $\delta_m$ both smaller than $1/m$, 
such that the following four properties hold:
\begin{equation}\delta_{m}=\delta_{f_{m}, a_{m}},\qquad \forall   m>1,
\end{equation}
\begin{equation}\label{sumbj}
\sum_{i=m}^\infty b_i(p)\leq  \delta_m(p),\qquad \forall   m>1,\forall p\in M
\end{equation}
\begin{equation}\label{almost isom}
(1-a_{m-1}) L_{g_m}(\gamma)\leq L_{\mathbb E}( f_m\circ\gamma)\leq L_{g_m}(\gamma)
,\qquad \forall \text{ curve }\gamma\subset M, \forall m>1, \text{ and}\end{equation}
\begin{equation}\label{dist-consec}
d_{\mathbb E}(f_{m-1}(p),  f_m(p))\leq b_{m-1}(p), \qquad \forall p\in M, \forall m>1  .
\end{equation}
Indeed, we already constructed $f_1$, $b_1$, and $\delta_1$. Assume that, for fixed $m$, 
$f_m$, $b_m$, and $\delta_m$ have been constructed. Let us construct 
$f_{m+1}$, $b_{m+1}$, and $\delta_{m+1}$.
Note that, since  $d_{g_m}\leq d_{g_{m+1}}$ and $f_m:(M,g_m)\to \mathbb E^k$ is $1$-Lipschitz, we have that $f_m:(M,g_{m+1})\to \mathbb E^k$ 
is $1$-Lipschitz as well. 
Applying Theorem \ref{Nash1} for $f_m$, $a_m$, and $b_m$, we get a 
$C^\infty$ $1$-Lipschitz embedding
$f_{m+1}:(M,g_{m+1})\to \mathbb E^k$ such that
\begin{equation*}
(1-a_m) L_{g_{m+1}}(\gamma)\leq L_{\mathbb E}( f_{m+1}\circ\gamma)\leq L_{g_{m+1}}(\gamma)
,\qquad \forall \text{ curve }\gamma\subset M, \end{equation*}
and
\begin{equation*}
d_{\mathbb E}(f_{m}(p),  f_{m+1}(p))\leq b_{m}(p), \qquad \forall p\in M  .
\end{equation*} 
Define $\delta_{m+1}=\delta_{f_{m+1}, a_{m+1}}$.
By induction, 
 we have that
$$\sum_{i=l}^m b_i < \delta_l,\qquad \forall l \text{ such that } 1\leq l\leq m .$$
Notice that  the above inequalities are strict.
Therefore 
 we can choose a    continuous function $b_{m+1}:M\to\R$ with $0<b_{m+1}<  \delta_{m+1}$ and such that
$$\sum_{i=l}^{ m+1} b_i < \delta_l,\qquad \forall l \text{ such that } 1\leq l\leq m+1 .$$
The construction of $\{f_m\}$, $\{b_m\}$, and $\{\delta_m\}$ is concluded.

We should notice that from \eqref{dist-consec} and \eqref{sumbj} we have that, if $m<j$,
 \begin{equation} \label{contNBHD0}
d_{\mathbb E}(f_{m}(p),  f_{j+1}(p))\leq \sum_{i=m}^j b_{i}(p)\leq \delta_m(p)\leq a_m=\dfrac{1}{m}.\end{equation}
In other words, for $j$ big enough,
 \begin{equation}
 f_{j}(M)\subset I_{\delta_m} (f_m(M))\subset \mathbb E^k.
 \label{contNBHD}
\end{equation}

After having constructed the sequence of approximating maps $f_m$, let us consider their limit. Notice that, since $d_{g_m}\leq d$, then the maps
$$f_m:(M,d)\to \mathbb E^k$$ 
are $1$-Lipschitz. 
By \eqref{contNBHD0}, the maps $f_m$ converge uniformly 
to a map
$$f:(M,d)\to \mathbb E^k,$$ 
which is obviously $1$-Lipschitz as well. Moreover, we have
 \begin{equation} \label{limcontNBHD2}
d_{\mathbb E}(f_{m}(p),  f (p)) \leq \delta_m(p)\leq a_m.\end{equation}

\subsubsection{The map {$f$} is a path isometry}
We will prove   that
 \begin{equation}
L_d(\gamma)\geq L_{\mathbb E}(f\circ\gamma), \qquad\forall \text{ curve }\gamma\subset M,
 \label{CL1}
\end{equation}
 and that
 \begin{equation}
L_d(\gamma)\leq L_{\mathbb E}(f\circ\gamma),\qquad \forall \text{ curve }\gamma\subset M.
 \label{CL2}
\end{equation}

The fact that \eqref{CL1} holds is obvious since $f$ is $1$-Lipschitz with respect to $d$.
For the proof of \eqref{CL2} we have to make use of the fact that $\delta_m$ have been constructed via the function $\delta$ of Lemma \ref{tubular}.
Observe that, taking limit in \eqref{contNBHD}, as $j\to\infty$, we have that, for all $m\in\N$,
 \begin{equation}
 f(M)\subset I_{\delta_m}  (f_m(M))\subset \mathbb E^k.
 \label{limcontNBHD}
\end{equation}
Let $I_m:=I_{\delta_m}  (f_m(M))$, and let $d_{I_m}$ be the path metric on it.

In order to prove \eqref{CL2}, take any curve $\gamma\subset M$ and take $p_0, p_1, \ldots, p_N\in\gamma$ consecutive points on the curve.
Fix   one of the indices $l\in\{1,\ldots,N\}$.
Consider the curve
$$\sigma_l:=[f_m(p_{l-1}),f(p_{l-1})]\cup f(\gamma|_{[p_{l-1},p_{l}]})\cup [f(p_{l}),      f_m(p_{l})],$$
where $[A,B]$, with $A,B\in\mathbb E^k$, is the Euclidean segment connecting $A$ and $B$.
By the \eqref{limcontNBHD2}, we have the 
containment
$$\sigma_l \subset I_m, \qquad \forall m\in\N.$$
In other words, the curve $\sigma_l$ connects the two points 
$f_m(p_{l-1})$ and $f_m(p_{l})$ inside the neighborhood $I_m$, so its length is greater than the path distance inside $I_m$ of such two points, i.e.,
$$d_{I_m}(f_m(p_{l-1}),f_m(p_{l}))\leq L_{\mathbb E}(\sigma_l).$$
Now, on the one hand, by the definition of $\sigma_l$ we have that
$$L_{\mathbb E}(\sigma_l)\leq 
\delta_m(p_{l-1})+L_{\mathbb E}( f\circ\gamma|_{[p_{l-1} ,p_{l}]})+\delta_m(p_{l})
\leq 2 a_m +L_{\mathbb E}( f\circ\gamma|_{[p_{l-1} ,p_{l}]}).$$
On the other hand, Lemma \ref{tubular} says that, since $\delta_m$  equals 
$\delta_{f_{m}, a_{m}}$, we have that
$$(1-a_m)d_{f_m(M)}( f_m(p_{l-1}),f_m(p_{l}))\leq d_{I_m}(f_m(p_{l-1}),f_m(p_{l})).$$
Therefore 
$$(1-a_m)d_{f_m(M)}( f_m(p_{l-1}),f_m(p_{l}))\leq 
2 a_m +L_{\mathbb E}( f\circ\gamma|_{[p_{l-1} ,p_{l}]}).$$
Since $f_m$ are $(1-a_{m-1})$-almost isometries (in the sense of \eqref{almost isom}),
we get
$$(1-a_m)(1-a_{m-1})d_{g_m}(  p_{l-1},p_{l})\leq 2 a_m +L_{\mathbb E}( f\circ\gamma|_{[p_{l-1} ,p_{l}]}).$$
Summing over $l$, we have that
$$(1-a_m)(1-a_{m-1}) \sum_{l=1}^Nd_{g_m}(  p_{l-1},p_{l})\leq 2 a_m N+L_{\mathbb E}( f\circ\gamma).$$
Now take the limit for $m\to\infty$. Since $a_m\to 0$,     (and note that $N$ is fixed), we get
$$  \sum_{l=1}^Nd (  p_{l-1},p_{l})\leq  L_{\mathbb E}( f\circ\gamma).$$
Finally, taking the supremum over all partitions of points $\{p_l\}$, we have that
$$L_d(\gamma)\leq L_{\mathbb E}(f\circ\gamma).$$

\subsubsection{The map {$f$} is an embedding} Assume by contradiction that there exists a point $q_0\in M$ and a sequence of points $q_k\in M$ with
$$f(q_k)\to f(q_0),\qquad \text{ but }\quad d(q_0,q_k)>\alpha, \forall k\in \N,$$
for some positive value $\alpha$. 
Since $d$ and $d_{g_{1}}$ give the same topology, there exists a $\beta>0$ such that
$$B_{d_{g_{1}}}(q_0,\beta)\subset B_d(q_0,\alpha).$$
Therefore, since the distances $d_{g_{ m}}$ are increasing, we can take $m$ large enough such that the following four inequalities hold:
\begin{equation}
d_{g_m}(q_0,q_k)\geq d_{g_{1}}(q_0,q_k)>\beta, \qquad \forall k\in\N
\end{equation}
\begin{equation}
1-a_m>\dfrac{1}{2},
\end{equation}
\begin{equation}
\delta_m<\dfrac{\beta}{16}, \text{ and}
\end{equation}
\begin{equation}
1-a_{m-1}>\dfrac{1}{2}.\end{equation}
Then, on the one hand,
\begin{eqnarray*}
d_{I_m} (f_m(q_k), f_m(q_0) )&\leq& d_{I_m} (f(q_k), f(q_0) )+  \delta_m (q_k)+  \delta_m (q_0)\\
&\leq& d_{I_m} (f(q_k), f(q_0) )+\dfrac{\beta}{8}.
\end{eqnarray*}
On the other hand,
\begin{eqnarray*}
d_{I_m} (f_m(q_k), f_m(q_0) )&\geq& (1-a_m)d_{f_m(M)}( f_m(q_k),f_m(q_{0}))\\
&\geq& (1-a_m)(1-a_{m-1})  d_{g_m}(  q_k,q_{0})
	\geq   \beta/4.
\end{eqnarray*}
So we get
$$d_{I_m} (f(q_k), f(q_0) ) \geq \dfrac{\beta}{4}-\dfrac{\beta}{8}=\dfrac{\beta}{8}>0,$$
which contradicts the fact that
$f(q_k)\to f(q_0)$, as $k\to\infty$.
\qed

\section{More on path isometric embeddings}\label{More}

\subsection{Sub-Riemannian geometries and the proof of Corollaries \ref{subRiemEmb} and \ref{subRiemEmb2}}
\label{subRiem Sect} 

\begin{definition}[The general definition of sub-Riemannian manifold]
 A (smooth) sub-Riemannian  structure on a manifold $M$ is a function $\rho :TM\to [0,\infty]$ obtained by the following construction: Let $E$ be a vector bundle over $M$ endowed with a scalar product  $\langle\cdot,\cdot\rangle$ and let
$$\sigma:E\to TM$$
be a morphism of vector bundles. For each $p\in M$ and $v,v'\in T_p M$, set 
$$\rho_p(v,v'):=\inf\{\langle u,u'\rangle: u,u'\in E_p, \sigma(u)=v, \sigma(u')=v'\}.$$
\end{definition}
Define $\rho_p(v):=\rho_p(v,v)$ and, given an absolutely continuous path $\gamma:[0,1]\to M$, define
$$L_\rho(\gamma):=\int_0^1 \sqrt{\rho_{\gamma(t)}(\dot\gamma(t))}d t.$$
The {\em sub-Riemannian distance associated to} $\rho$ is defined as, for any $p$ and $q$ in $M$,  
$$d_{CC}(p,q)=\inf\left\{ L_\rho(\gamma)\;\Big{|}\; \gamma \text{ absolutely continuous path } \gamma(0)=p, \gamma(1)=q\right\}.$$
The only extra assumption on $\rho$ is that the distance $d_{CC}$ is finite and induces the manifold topology.

\proof[Proof of Corollary \ref{subRiemEmb}]
We show now  that each sub-Riemannian distance can be obtained as a limit of increasing Riemannian distances. 
The proof is easy and well-known in the case when  $E$ is in fact a
sub-bundle of the tangent bundle. Here we give the  proof in the general case.

Let  $\rho :TM\to [0,\infty]$  be the function defining the sub-Riemannian structure.
Notice that $\rho(v)=0$ only if $v=0$. So one can take some Riemannian tensor $g_{1}$ with the property that 
$g_{1}\leq \rho.$

Then, by recurrence, for each $m\in \N$, we consider $g_m$ to be a (smooth)  Riemannian tensor with the property that, at any point $p\in M$, 
$$\max\{(g_{m-1})_p(v,w),\min\{ (1-2^{-m})\rho_p(v,w),m(g_{1})_p(v,w)\}\}\leq (g_{m})_p(v,w)\leq \rho_p(v,w) .$$
Obviously  we have that
$$g_{1}\leq g_m\leq g_{m+1}\leq\rho.$$
Then, for any absolutely continuous path $\gamma$, we have that
$$L_{g_m}(\gamma)\leq L_{\rho}(\gamma) .$$
Thus, for any $p$ and $q$ in $M$, 
$$d_{g_m}(p,q)\leq d_{CC}(p,q) ,$$
and therefore
$$\lim_{m\to\infty}d_{g_m}(p,q)\leq d_{CC}(p,q) .$$
Assume, by contradiction, that, for some $p$ and $q$ in $M$,  we have that
$$\lim_{m\to\infty}d_{g_m}(p,q)< d_{CC}(p,q) .$$
Then there are curves $\gamma_m$ from $p$ to $q$ such that
$$\lim_{m\to\infty}L_{g_m}(\gamma_m)< d_{CC}(p,q) .$$
Since 
$$L_{g_{1} }(\gamma_m)\leq L_{g_m}(\gamma_m),$$
we get a bound on the lengths $L_{g_{1} }(\gamma_m)$. Therefore, by an Ascoli-Arzel\`a argument, $\gamma_m$  converges to a curve $\gamma$  from $p$ to $q$.
We may assume that $\gamma$ is parametrized by arc length with respect to the distance of $g_1$.
Now, either $L_\rho( \gamma)$ is infinite or it is finite.
Namely, either there is a positive-measure set $A\subset[0,L_{g_1}(\gamma)]$ such that
$$\rho_{\gamma(t)}(\dot\gamma(t))=\infty, \qquad\forall t\in A,$$
or, for almost every $t\in [0,L_{g_1}(\gamma)]$, the value
$\rho_{\gamma(t)}(\dot\gamma(t))$ is finite.

In the first case,
for all $t\in A$, 
$$(g_m)_{\gamma(t)} (\dot\gamma(t))\geq  m \cdot (g_{1} )_{\gamma(t)} (\dot\gamma(t)).$$
From this we have that
$$L_{g_m}(\gamma)\geq m L_{g_ {1}}(\gamma|_A)\to \infty, \quad \text{ as } m\to \infty.$$
We get a contradiction since by assumption  $d_{CC}(p,q)<\infty$.

In the second case,
for almost all $t$, for $m$ big enough, 
$$(1-2^{-m})  \rho_{\gamma(t)} (\dot\gamma(t))\leq    (g_m)_{\gamma(t)} (\dot\gamma(t))\leq     \rho_{\gamma(t)} (\dot\gamma(t)).$$
From this we have that
$$L_{g_m}(\gamma)\to  L_{\rho}(\gamma), \quad \text{ as } m\to \infty.$$
We get a contradiction since we have that 
   $d_{CC}(p,q)\leq L_{\rho}(\gamma)$.\qed

%
%

\proof[Proof of Corollary \ref{subRiemEmb2}]
Corollary \ref{subRiemEmb2} is not a direct consequence of the claim of Theorem \ref{path-isom-embed}. However, the proof is the same.
Indeed, in the proof of the theorem we started with the embedding
$$f_1:(M,g_1)\to \mathbb E^k$$
with $k=2 n+1$, which was given by Theorem \ref{Nash2}.
If instead, as assumed in Corollary \ref{subRiemEmb2}, we already have an embedding
$$f:(M, d_{\rm Riem})\to \mathbb E^k$$
with $k\geq n+1$, then we can consider a sequence of increasing Riemannian distances starting with $d_{g_1}=d_{\rm Riem}$ and converging point-wise to $d$.
At each stage, each $1$-Lipschitz embedding can be stretched as in Theorem \ref{path-isom-embed},  since in Theorem \ref{Nash1} we only need the codimension to be greater than $1$, i.e., 
$k\geq n+1$.
\qed

\begin{corollary}\label{HinR4}
 Let $(\mathbb H,d_{CC})$ be the Heisenberg group endowed with the  sub-Riemannian distance with the first layer as horizontal distribution.
Then we have that 
there exists a subset $\Sigma$ of $\R^4$, such that, if $d_\Sigma$ is the path metric induced by the Euclidean length of  $\R^4$,
then $(\mathbb H,d_{CC})$  is isometric to $(\Sigma,d_\Sigma)$.
\end{corollary}
\proof The statement is a direct consequence of Corollary \ref{subRiemEmb2} and Proposition \ref{path-isom-isom}. We make use of the fact that 
 the inverse of the stereographic projection, which   maps  $\R^3$ to $\mathbb S^3\subset \R^4$, gives a globally
 Lipschitz embedding  of the
 Riemannian left-invariant Heisenberg group into the Euclidean space  $\mathbb E^4$.
\qed
\begin{remark} 
A similar reasoning can be applied to the Grushin plane. The reader can be referred to \cite{Monti-Morbidelli} for an introduction to the geometry of the Grushin plane.
Thus, another 
 consequence of Corollary \ref{subRiemEmb2} and Proposition \ref{path-isom-isom}
is the following fact. 
The Grushin plane $\mathbb P$ can be realized as a subset of $\R^3$ with the induced path distance. The reason is again that  the inverse of the stereographic projection 
gives a Lipschitz embedding
of $\mathbb P$ into 
$\mathbb E^3$. 
 \end{remark}

\subsection{Isometries, intrinsic isometries, and path isometries}\label{Isometries}
\label{intrinsic isometry section}
This section is devoted to the equivalence of the various  notions of path isometric embeddings and of intrinsic isometric embeddings. 
\begin{proposition}\label{path-isom-isom}
 Let $f:(X,d_X)\to (Y,d_Y)$ be a map between proper geodesic metric spaces.
Then $f$ is a path isometric embedding if and only if the space $f(X)$ endowed with the path distance $d_{f(X)}$ induced by $d_Y$ is isometric to $(X,d_X)$ via $f$
and 
 the topology induced by $d_{f(X)}$ coincides with the topology of $f(X)$ as a topological subspace of $Y$.
\end{proposition}
\proof
Let us denote by $\tau_X$ and $\tau_Y$ the topology of $(X,d_X)$ and $ (Y,d_Y)$, respectively.
Let $\tau_{d_{f(X)}}$ be the topology on $f(X)$ induced by the path distance $d_{f(X)}$.
We shall write $A\simeq B$ to say that $A$ is homeomorphic to $B$.

$\Leftarrow]$ If $f:(X,d_X)\to (f(X),d_{f(X)})$ is an isometry, then it preserves the length of paths. Since the length structures on $f(X)$ and $Y$ coincide, then $f:(X,d_X)\to (Y,d_Y)$ is a path isometry, cf. \cite[Proposition 2.3.12]{Burago}.
Moreover, since $f:(X,d_X)\to (f(X),d_{f(X)})$ is an isometry, then $(X,\tau_X)\simeq (f(X),\tau_{d_{f(X)}}  )$.
If, by assumption  $(f(X),\tau_Y) \simeq (f(X),\tau_{d_{f(X)}}  )$, we have that  $(f(X),\tau_Y) \simeq(X,\tau_X)$, i.e., $f$ is an embedding.

$\Rightarrow]$ If $f$ is an embedding, we have, by definition, that $(f(X),\tau_Y) \simeq(X,\tau_X)$.
Moreover, since $f$ has a continuous inverse on $f(X)$,  there is a one-to-one correspondence
between curves in $X$ and curves in $f(X)$.
If $f$ is a path isometry, then such a correspondence preserves length. Since both $d_X$ and $d_Y$ are length spaces, we have that  
$$d_X(x,y)=d_{f(X)}(f(x),f(y)), \qquad x,y\in X,$$
i.e., $f:(X,d_X)\to (f(X),d_{f(X)})$ is an isometry.

We also have as a consequence that $(X,\tau_X)\simeq (f(X),\tau_{d_{f(X)}}  )$. 
 We conclude that
 $(f(X),\tau_Y) \simeq (f(X),\tau_{d_{f(X)}}  ).$
\qed

We recall now the definition of intrinsic isometry. The aim is to relate our work with the one of Petrunin \cite{Petrunin}. 
Let $f:X\to Y$ be a map between  length spaces. 
Given two points $p,q\in {X}$, a sequence of points $p=x_0,x_1,\dots,x_N=q$ in $X$ is called an {\em $\eps$-chain from $p$ to $q$} if $d(x_{i-1},x_i)\le\eps$ for all $i=1,\ldots,N$.
Set
$${\rm pull}_{f,\eps}(p,q)
=
\inf\left\{\sum_{i=1}^N d(f(x_{i-1}),f(x_{i}))\right\},$$
where the infimum is taken along all  $\eps$-chains $\{x_i\}_{i=0}^N$ from $p$ to $q$.
The limit
$${\rm pull}_{f}(p,q):=\lim_{\eps\to0}{\rm pull}_{f,\eps}(p,q)$$
defines a (possibly infinite) pre-metric.

A map $f:{X}\to {Y}$  is called an \emph{intrinsic isometry} if 
$$d_X(p,q)={\rm pull}_{f}(p,q)$$
for any $p,q\in {X}$.

\begin{proposition}\label{intrinsic prop}
 A path isometric embedding $f:X\to Y$ between proper geodesic spaces is an intrinsic isometry.
\end{proposition}
\proof
Take $p$ and $q\in X$. Let $\gamma$ be a geodesic from $p$ to $q$. Fix $\eps>0$. Let $t_0<t_1<\ldots<t_N$ be such that 
$$\gamma(t_0)=p,\qquad \gamma(t_N)=q,$$
and
$$\{\gamma(t_j)\}_{j=0}^N \text{ is an } \eps\text{-chain}.$$
Then, using that $f$ is a path isometry, we have that
\begin{eqnarray*}
 {\rm pull}_{f,\eps}(p,q)&\leq& \sum_{i=1}^N d(f(\gamma(t_{i-1})),f(\gamma(t_{i})))\\
&\leq&\sum_{i=1}^N L_Y(f\circ\gamma|_{[t_{i-1},t_{i}]}  )\\
&=&\sum_{i=1}^N L_X(\gamma|_{[t_{i-1},t_{i}]}  ) 
= L_X(\gamma)  
=d(p,q).
\end{eqnarray*}
To prove the other inequality, assume by contradiction that there is some $\alpha>0$ and there is some $\eps_0>0$ such that, for all $\eps\in(0,\eps_0)$,
we have that
$$ {\rm pull}_{f,\eps}(p,q)\leq d(p,q)-\alpha.$$
Thus, for each such an $\eps$ there exists an $\eps$-chain $\{x_i^{(\eps)} \}_{i=0}^N$ from $p$ to $q$ with the property that
$$\sum_{i=1}^N d(f(x_{i-1}^{(\eps)}),f(x_{i}^{(\eps)}))\leq d(p,q)-\alpha/2.$$
Consider a curve $\sigma_\eps$ in $Y$ passing  through the points  $f(x_0^{(\eps)}), f(x_1^{(\eps)}),\dots,f(x_N^{(\eps)})$ and forming a geodesic between 
$f(x_{i-1}^{(\eps)}) $ and $  f(x_{i}^{(\eps)})$.
Therefore we have that
$$L_Y(\sigma_\eps) \leq d(p,q)-\alpha/2.$$
From such a bound on the length, from the fact that $\sigma_\eps$ starts at the fixed point $f(p)$, and from the fact that $Y$ is locally compact, we have that there exists a limit curve $\sigma$, as $\eps\to 0$, with the property that 
$$L_Y(\sigma) \leq d(p,q)-\alpha/2.$$
Since $\{f(x_{i}^{(\eps)}) \}_{i=0}^N$ are finer and finer on $\sigma_\eps$, as $\eps\to0$, we have 
$\sigma \subset f(X)$.
Since $f$ is a homeomorphism between $X$ and $f(X)$, we have the existence of a curve $\gamma$ from $p$ to $q$ with the property that
$$f\circ\gamma=\sigma.$$
We arrive at a contradiction since
\begin{eqnarray*}
 d(p,q)&\leq&L_X(\gamma) 
 =L_Y(\sigma)  
 \leq d(p,q)-\alpha/2.
\end{eqnarray*}
\qed

\subsection{Metric spaces that are not path isometrically embeddable}

\proof[Proof of Proposition \ref{Finsler}.]
 We prove that the norm $\norm{\cdot}$ at a point comes from a scalar product by showing that it is the pull back norm of an Euclidean norm via a linear map.
Roughly speaking, we would like to claim the following. Assume that $f$  is differentiable at $p$.
Since $f$ is a path isometry,  it  sends infinitesimal metric balls at $p$ in $(M,\norm{\cdot})$ to infinitesimal metric balls at $f(p)$ in $(f(M),d_{\mathbb E})$.
However,    infinitesimal balls at $f(p)$ are spheres and, $df_p$ being linear, 
infinitesimal balls at $p$ would be ellipsoids.

Consider an open set $U\subset \R^n$ and a smooth coordinate chart $\phi:U\to M$.
Notice that  $f:(U,d_ \mathbb E)\to (M,\norm{\cdot})$ is locally  Lipschitz.

 If 
$f:(M,\norm{\cdot})\to \mathbb E^k,$
 is a path isometry, then it is a $1$-Lipschitz map. 
 Hence $F:=f\circ \phi$ is locally a Lipschitz map between Euclidean open sets.
 According to Rademacher's Theorem, for almost all $q\in U$, 
 the differential $dF_q$ exists and the map $v\mapsto dF_q v$   is linear.
  We fix a dense and countable set of directions $\mathcal D \subset \R^n$. 
  Hence, considering Lebesgue points of the measurable functions  $q\to    dF_q v$, we obtain that, for almost all $q\in U$ and all directions $v\in \mathcal D$, the differential $dF_q$ exists  and is linear and
  \begin{equation}\label{ratio}
\lim_{\eps\to 0} \dfrac{L_{\mathbb E} ( F(q+t v)|_{t\in[0,\eps]} )}{\norm{F(q+\eps v)- F(q)}_{\mathbb E}}
=
\lim_{\eps\to 0} \dfrac{ \int_0^\eps  \norm{ \frac{d}{dt} F(q+t v)}_{\mathbb E}  \,dt}{\norm{\int_0^\eps  \frac{d}{dt} F(q+t v)\,dt }_{\mathbb E}  }
=1.  
  \end{equation}

%
%
%


 Since
 $\norm{\cdot}$ is smooth and the curve
 $t\mapsto \phi(q+tv)$ is smooth,
 we have
 $$ \norm{d\phi_q v}=\lim_{\eps\to 0}\frac{1}{\eps}L_{\norm{\cdot}} ( \phi(q+tv)|_{t\in[0,\eps]} ).$$
Since $f$ is a path isometry, the latter equals 
$$\lim_{\eps\to 0}\frac{1}{\eps}L_{\mathbb E} ( (f\circ \phi)(q+tv)|_{t\in[0,\eps]} ).$$
 If $q$ is one of the above points where $F=f\circ \phi$ is differentiable and \eqref{ratio} holds with $v\in \mathcal D$, then
\begin{eqnarray*}\lim_{\eps\to 0}\frac{1}{\eps}L_{\mathbb E} ( F(q+tv)|_{t\in[0,\eps]} )  
  &=&
  \lim_{\eps\to 0}\frac{1}{\eps}
  \norm{F(q+\eps v)- F(q)}_{\mathbb E}
  \\
&=&\norm{(dF_q)(v)}_{\mathbb E}.\end{eqnarray*}
Since the set of directions $\mathcal D$ is dense, we get
$$\norm{d\phi_q v}=\norm{(dF_q)(v)}_{\mathbb E}, \qquad \forall v\in T_q\R^n.$$
In other words, $\norm{\cdot}$ at $q$ is the pull back norm via $dF_q$ of the Euclidean norm $\norm{\cdot}_{\mathbb E}$.
Since $dF_q$ is linear, the norm $\norm{\cdot}$ at $q$ comes from  a scalar product.
Since we can consider a sequence of points  $\phi(q)$ tending to $p$, we also have the same result for the generic $p$, by continuity of the Finsler structure.
\qed

\section{Lipschitz embeddings for finite dimensional metric spaces}\label{Lipschitz}
\subsection{Preliminaries}
This section is a preparation to the proof of the Embedding Theorem \ref{Lip-embed}. 
To fix some notation, we recall the notion of general position.  
A set $\lbrace \mathbf{x}_0, \ldots , \mathbf{x}_k\rbrace$ of points of $\R^N$  is said to be {\em geometrically independent}, or  {\em affinely independent}, if the equations
$$ \sum_{j=1}^k a_j\mathbf{x}_j=\mathbf{0} \qquad {\rm and} \qquad \sum_{j=1}^k a_j=0$$
hold only if each $a_j =0$.
In the language of ordinary linear algebra, this is just the definition of linear independence for the set of vectors $\mathbf{x_1} - \mathbf{x}_0, \ldots, \mathbf{x}_k - \mathbf{x}_0$ of the vector space $\R^N$.
So $\R^N$ contains no more than $N + 1$ geometrically independent points.


   A set $A$ of points of $\R^N$ is said to be in {\em general position in} $\R^N$ if every subset of $A$
   containing $N + 1$ or fewer points
     is geometrically independent.
Observe that, given a finite set $\{\mathbf{x}_1, \ldots,\mathbf{x}_n\}$ of points of $\R^N$ 
and given $\delta > 0$, there exists a set $\{\mathbf{y}_1, \ldots, \mathbf{y}_n\}$ of points of $\R^N$ in general position in $\R^N$, such that $|\mathbf{x}_j - \mathbf{y}_j| < \delta$ for all $j$.
\\



{\prop\label{projection} Suppose $K$ is a compact subset of $\R^n$ of packing dimension $k$. If $n > 2k + 1$, then there
is a full measure subset $A$ of the unit
sphere $\mathbb S^{n-1}$ such that if $v$ is an element of $A$, and
$$\pi_v : \R^n \To \R^{n-1}$$
is the orthogonal projection along $v$, then  the restriction of
$\pi_v$ to $K$ is a (Lipschitz) homeomorphism.}

\proof The proof is based on the fact that every pair of distinct points
in $K$ determines a line in $\R^n$, and hence an element of projective
space $\R P^{n-1} = \mathbb S^{n-1}/\{\pm 1\}$.  
Recall that the Hausdorff dimension of  $K \times K $ is bounded by twice the 
packing dimension of $K$.
The map $K \times K \setminus {\rm Diag}(K\times K) \To \R P^{n-1}$
is locally Lipschitz. Thus its image has Hausdorff dimension $\leq 2k$.
The complement in $\R P^{n-1}$ gives the set $A$.
\qed

\rem We can iterate the proposition to conclude that, if $K$ is a compact $k$-dimensional  subset of $\R^n$, we can find a (full-measure) set of orthogonal projections $\tilde\pi:\R^{n}\To \R^m$, as soon as $n>m = 2k + 1$, that are homeomorphisms when restricted to $K$.

\rem Since $A$ has full measure, it is    dense. Thus, given any projection, it is possible to find a `good' projection as close as we want. 
\\

The core of the proof in the theorem of Menger and  N\`obeling is the construction  of embeddings that are close to being injective. One  uses the analytic geometry of $\R^N$ discussed earlier. We present now the relative version for the Lipschitz case.

{\lem\label{goodg} If $(X,d)$ is a compact metric space of topological dimension $m$, then, for all $N\geq 2m+1$, there exists a Lipschitz map 
arbitrarily close to being injective with range into the Euclidean space of dimension 
$N$, i.e., for any
fixed 
$\eps>0$ there exists $g\in {\rm Lip}(X; \R^N)$ such that
$$ g(x_1)=g(x_2)\Longrightarrow d(x_1,x_2)<\eps.$$}

\proof By the definition of topological dimension, we have that we can cover $X$ by finitely many open sets $\{U_1, \ldots, U_n\}$ such that
\begin{enumerate}
\item $\diam U_j <\eps$ in $X$,
\item $\{U_1, \ldots, U_n\}$ has order $\leq m + 1$.
\end{enumerate}
The second requirement means that no point of $X$ lies in more than $m+1$ elements of the cover.

Let $\phi_j$ be a Lipschitz partition of unity dominated by $\{U_j\}$, cf.~\cite{Luukkainen-Vaisala}. 
For each $j$, choose a point  
$\mathbf{z}_j\in \R^N$ 
such that the set $\{\mathbf{z}_1, \ldots, \mathbf{z}_n\}$ is in general position in $\R^N$. Finally, 
define $g : X \To \R^N$ by the equation
$$ g(x)=\sum_{j=1}^n\phi_j(x)\mathbf{z}_j.$$
We assert that $g$ is the desired function.

At every point $x$, locally $g(x)$ is a sum of finitely many Lipschitz maps, thus is Lipschitz. 

We shall prove that if $x_1$, $x_2\in X$ and $g(x_1) = g(x_2)$, then $x_1$ and $x_2$ belong to one of the open sets $U_j$, 
so that necessarily $d(x_1, x_2)<\eps$ (since $\diam U_j < \eps$). 

So suppose $g(x_1) = g(x_2)$. Then
$$\sum_{j=1}^n\left[ \phi_j(x_1)-\phi_j(x_2)\right] \mathbf{z}_j=0.$$
Because the covering $\{ U_j \}_{j=1}^n$ has order at most $m + 1$, at most $m + 1$ of the numbers $\{\phi_j(x_1)\}_{j=1}^n$ are nonzero, and at most $m + 1$ of the numbers $\{\phi_j(x_2)\}_{j=1}^n$ are nonzero. Thus, the sum $$\sum_{j=1}^n \left[ \phi_j(x_1)-\phi_j(x_2)\right] \mathbf{z}_j=0$$ has at most $2m + 2$ nonzero summands. Note that the sum of the coefficients vanishes because
$$\sum_{j=1}^n\left[ \phi_j(x_1)-\phi_j(x_2)\right]= 1-1 = 0.$$
The points $\mathbf{z}_j$, are in general position in $\R^N$, so that any subset of them having $N + 1$ or fewer elements is geometrically independent. And by hypothesis $N + 1 = 2m + 2$.  Therefore, we conclude that
$$\phi_j(x_1)-\phi_j(x_2) = 0$$
for all $j$.
Now $\phi_j(x_1)> 0$ for some $j$, so that $x_1\in U_j$. Since $\phi_j(x_1)-\phi_j(x_2)=0$, we have that $x_2\in U_j$ as well, as asserted.
\qed

\subsection{The proof of the Embedding Theorem \ref{Lip-embed}}
 Let $X$ be a compact metric space of finite packing dimension. Let $k$ be the  packing dimension of $X$.
    Let $m$ be the  topological dimension of $X$. Hence, $m\leq k$.
   Set $N := 2k + 1$. 

Consider the space ${\rm Lip}(X; \R^N)$, i.e., the space of all the Lipschitz maps from $X$ to $\R^N$.
It is non-empty, since the constant functions are there.
It is complete in the following metric: 
$$\norm{f}_{Lip}:=\norm{f}_\infty+\sup\left\lbrace \dfrac{|f(x)-f(y)|}{d(x,y)}\,:\,x,y\in X, x\neq y\right\rbrace .$$

Let $d$ be the metric of the space $X$; because $X$ is compact, $d$ is bounded. Given a map $f : X \To \R^N$, let us define
$$\Delta(f) := \sup\{\diam f^{-1}(\mathbf{z}) :  \mathbf{z}\in \R^N\},$$
i.e., the fibers of $f$ have diameter smaller than $\Delta(f)$.
So the number $\Delta(f)$  measures how far $f$ is far from being injective; if $\Delta(f)= 0$, then in fact $f$ is injective.

Now, given $\varepsilon> 0$, define $\mathcal{U}_\varepsilon$ to be the set of all those Lipschitz maps $f: X\To \R^N$ for which $\Delta(f)< \eps$. In Lemma \ref{open} and in Lemma \ref{dense} we shall show that $\mathcal{U}_\eps$ is both open and dense in ${\rm Lip}(X; \R^N)$, respectively. So it follows from the Baire Category  Theorem that the intersection
$$ \bigcap _{n\in\N}\mathcal{U}_{\frac{1}{n}}$$
is dense in  ${\rm Lip}(X; \R^N)$ and is in particular non-empty.
If $f$ is an element of this intersection, then $\Delta(f)< 1/n$ for every $n$. Therefore, $\Delta(f)=0$ and $f$ is injective. Because $X$ is compact, $f$ is an embedding. Thus, modulo   Lemma  \ref{open} and Lemma \ref{dense}, the theorem is proved. \qed

\begin{lemma}\label{open} $\mathcal{U}_\eps$ is open in ${\rm Lip}(X; \R^N)$.\end{lemma}
\proof
Given an element $f\in\mathcal{U}_\eps$, we wish to find a ball at $f$ of some radius $\delta$ that is contained in $\mathcal{U}_\eps$. First choose a number $b$ such that $\Delta(f)< b<\eps$. 
Let $A$ be the following subset
$$A = \left\lbrace (x, y) \in X \times X \, : \, d(x, y) \geq  b\right\rbrace .$$
Now $A$ is closed in $X \times X$ and therefore compact.

Note that if $f(x) = f(y)$, then $d(x , y)$ must be less than $b$. 
It follows that the function $|f(x) - f(y)|$ is positive on $A$. Since $A$ is compact, the function has a positive minimum on $A$.  Let
$$\delta := \frac{1}{2}\min\left\lbrace |f(x) - f(y)| \, : \, x,y\in A\right\rbrace.$$
We assert that this value of $\delta$ will suffice.

Suppose that $g$ is a map such that $\norm{f-g}_{Lip} < \delta$, so in particular $\norm{f-g}_{\infty} < \delta$. 
If $(x, y) \in A$, then $|f(x) - f(y)| > 2\delta$ by definition of $\delta$.
Since $g(x)$ and $g(y)$ are within $\delta$ of $f(x)$ and $f(y)$, respectively, we  must have that $|g(x) - g(y)| > 0$. 
Hence the function $|g(x) - g(y)|$ is positive on $A$. As a  result, if $x$ and $y$ are two points such that $g(x) = g(y)$, then necessarily $d(x, y) < b$. We conclude that $\Delta(g)\leq b<\eps$, as desired.
\qed
\begin{lemma}\label{dense} $\mathcal{U}_\eps$ is dense in ${\rm Lip}(X; \R^N)$.\end{lemma} 
\proof
This is the more substantial  part of the proof. We shall use the preliminaries presented in the previous subsection. Let $f\in {\rm Lip}(X; \R^N)$. Given $\delta > 0$
, we wish to find a function $F\in {\rm Lip}(X; \R^N)$
such that $F\in\mathcal{U}_\eps$ and $\norm{f-F}_{Lip}<\delta$.

Since the topological dimension $m$ of $X$ is at most $k$, we can apply Lemma \ref{goodg}. Take $g\in {\rm Lip}(X; \R^N)$
such that  
if $ g(x_1)=g(x_2)$ then $d(x_1,x_2)<\eps/2$. 

Consider $\Phi := (f, g):X\To \R^{2N}$.
Clearly, $\Phi$ is Lipschitz. Thus, $\Phi(X)$ has packing dimension no more than $k$.

Since $2N>N=2k+1$, we can use Proposition \ref{projection} (and the remarks afterwards) to project   the compact set $K=\Phi(X)$  
 from $\R^{2N}$ to $ \R^N$. Namely, there are orthogonal projections that are injective on $K$ and are arbitrarily close to the projection in the first $N$-dimensional component.
Explicitly, for any $\beta>0$, there exists an orthogonal projection  $\tilde\pi:\R^{2N}\To \R^N$ such that 
the restriction of
$\tilde\pi$ to $K$ is a (Lipschitz) homeomorphism and,
if $\pi:\R^{2N}=\R^N\times\R^N\To \R^N$ is given by $\pi(\mathbf{x},\mathbf{y})=\mathbf{x}$,
then $$ \norm{\tilde\pi-\pi}< \beta.$$ We are using here the operator norm. We will say later how small $\beta$ has to be  in terms of the data $(f,g,\delta)$. 

Set $F:=\tilde\pi\circ \Phi$. We shall prove first that $F\in\mathcal{U}_\eps$ and then $\norm{f-F}_{Lip}<\delta$.

Suppose $x_1$, $x_2$ are in the same fiber of $F$, i.e., $F(x_1)=F(x_2)$. So from the definition of $F$,
$(\tilde\pi\circ \Phi)(x_1)=(\tilde\pi\circ \Phi)(x_2).$
Since $\tilde\pi$ is a homeomorphism on $K=\Phi(X)$, we have that 
$\Phi(x_1)=\Phi(x_2).$
From the definition of $\Phi$, we have that 
$$\left( f(x_1), g(x_1)\right)=\left(f(x_2), g(x_2)\right).$$
In particular, $g(x_1) = g(x_2)$. So, by the property of $g$, we have that $d(x_1,x_2)<\eps/2$.
Therefore, $F\in\mathcal{U}_\eps$.

Let us prove now that $F$ is $\delta$-close to $f$. Let us write explicitly the difference
\begin{eqnarray*}
 F(x)-f(x)		& = & (\tilde\pi\circ \Phi)(x)-f(x) \\
			& = & \tilde\pi\left(f(x),g(x)\right)-\pi\left(f(x),g(x)\right)
				=(\tilde\pi-\pi)\left(f(x),g(x)\right).
\end {eqnarray*}
Bound the sup norm by
\begin{eqnarray*}
|F(x)-f(x)|		&\leq& \norm{\tilde\pi-\pi}|\left(f(x),g(x)\right)|\\
			&\leq& \norm{\tilde\pi-\pi} \sqrt{\norm{f}^2_\infty+\norm{g}^2_\infty}  
			\leq \beta\sqrt{\norm{f}^2_{Lip}+\norm{g}^2_{Lip}}.
\end {eqnarray*}
For the Lipschitz part of the norm, remember that the projections are linear. Therefore
\begin{eqnarray*}
\dfrac{|F(x)-f(x)-(F(y)-f(y))|}{|d(x,y)|}&\leq&\dfrac{| (\tilde\pi-\pi)\left(f(x),g(x)\right)-(\tilde\pi-\pi)\left(f(y),g(y)\right)|}{d(x,y)}\\
&\leq&\dfrac{| (\tilde\pi-\pi)\left(f(x)-f(y),g(x)-g(y)\right)|}{d(x,y)}\\
&\leq&\norm{\tilde\pi-\pi}\dfrac{|\left(f(x)-f(y),g(x)-g(y)\right)|}{d(x,y)}\\
			&\leq& \norm{\tilde\pi-\pi} \sqrt{\norm{f}^2_{Lip} +\norm{g}^2_{Lip}}\\
			&\leq& \beta\sqrt{\norm{f}^2_{Lip}+\norm{g}^2_{Lip}}.
\end {eqnarray*}
So choose $\beta$ such that $\beta\sqrt{\norm{f}^2_{Lip}+\norm{g}^2_{Lip}}<\delta/2$.
\qed

\section{Acknowledgment}
The author is indebted to the stimulating remarks, helpful advice, and useful suggestions of B. Kleiner and A. Petrunin.

 \bibliography{general_bibliography}
\bibliographystyle{amsalpha}
 

  \vskip 1in

\parbox{3.5in}{Enrico Le Donne:\\
~\\
enrico.ledonne@math.ethz.ch\\
enrico.ledonne@math.u-psud.fr\\ 
enrico.ledonne@msri.org
}

\end{document}